\documentstyle[11pt]{article}
\begin{document}
\newtheorem{prop}{Proposition}[section]
\newtheorem{Def}{Definition}[section]
\newtheorem{theorem}{Theorem}[section]
\newtheorem{lemma}{Lemma}[section]
\newtheorem{Cor}{Corollary}[section]

\title{\bf Bounds in time for the Klein - Gordon - Schr\"odinger and the Zakharov system}
\author{{\bf Axel Gr\"unrock}\\{\bf Hartmut Pecher}\\Fachbereich Mathematik und Naturwissenschaften\\ Bergische Universit\"at Wuppertal\\ Gau{\ss}str. 20 \\ D-42097 Wuppertal\\ Germany\\ e-mail Axel.Gruenrock@math.uni-wuppertal.de\\Hartmut.Pecher@math.uni-wuppertal.de}
\date{}
\maketitle

\begin{abstract}
It is shown that the spatial Sobolev norms of regular global solutions of the (2+1), (3+1) and (4+1)-dimensional Klein-Gordon-Schr\"odinger system and the (2+1) and (3+1)-dimensional Zakharov system grow at most polynomially with a bound depending on the regularity class of the data. The proof uses the Fourier restriction norm method.
\end{abstract}

\normalsize
\setcounter{section}{-1}
\section{Introduction}
Consider the Cauchy problem for the Klein - Gordon - Schr\"odinger system (KGS) with Yukawa coupling in space dimension $2 \le n \le 4 $ : 
\begin{eqnarray}
i\psi_t + \Delta \psi & = & -\phi\psi \nonumber\\
\phi_{tt} + (-\Delta + 1)\phi & = & |\psi|^2 \label{2.1} \\
\psi(0) \; = \; \psi_0 \quad , \quad \phi(0) & = & \phi_0 \quad , \quad \phi_t(0)\; = \;\phi_1 \nonumber
\end{eqnarray}
where $\psi$ is the complex-valued nucleon field and $\phi$ the real-valued meson field.

Local existence and uniqueness holds for data $(\psi_0,\phi_0,\phi_1) \in H^{s} \times H^{m} \times H^{m-1},$ if $ s \ge 1 $ , $ m > 0 $ and $ s-1<m<s+2 $ , and the life span depends only on $ \| \psi_0\|_{H^{1}} $ , $ \|\phi_0\|_{H^{\delta}} $ and $ \|\phi_1\|_{H^{\delta -1}}$ , where $ \delta = 0 $ , if $ n\le 3$ ,  and $\delta > 0$ arbitrarily small, if $ n=4 $ (for details see Prop. \ref{Proposition 2.1} below). Cf. also \cite{P}, Thm. 2.1.

If we assume data with finite energy, i.e. $s,m\ge 1$, the conservation laws imply that this local solution can be extended globally in time, if $n=2$ or $n=3$, and, provided $\|\psi_0\|_{L^2}$ is small, if $n=4$.

So the question arises whether a non-trivial bound in time for the Sobolev norms of the solution can be given. We show in Theorem \ref{Theorem 2.1} below that for all $ s >> 1 $ a polynomial bound can be given, and especially for $s\le m \le s+1$ we have the estimate
$$ \|\psi(t)\|_{H^{s}} + \|\phi(t)\|_{H^{m}} + \|\phi_t(t)\|_{H^{m-1}} \le c (1+t)^{s-1} $$
To prove such a result we use the Fourier restriction norm method using bilinear estimates in order to control the nonlinear terms and to get an estimate of the form
$$  f(t) \le f(\tau) + c f(\tau)^{1-\delta} $$
with $ 0 < \delta \le 1 $ on any interval $ t \in [\tau,\tau+t_0] $ for a fixed $t_0$. Here $f(t)$ is the sum of the involved Sobolev norms of $\psi$ , $ \phi $ and $\phi_t$.
Bourgain \cite{B} observed that such an inequality implies the bound (cf. Lemma \ref{Lemma 2.2})
$$ f(t) \le c_1 + c_2 t^{\frac{1}{\delta}} $$
and used it to prove polynomial bounds for Sobolev norms of solutions to nonlinear Schr\"odinger and wave equations \cite{B}.

Staffilani \cite{St1},\cite{St2} improved these results in the case of the nonlinear Schr\"odin-\\ger equation and considered moreover KdV-type equations.

Colliander, Delort, Kenig and Staffilani \cite{CDKS} improved the results for the (2+1)-dimensional Schr\"odinger equation $iu_t + \Delta u \pm |u|^2 u = 0$ further and were able to show the bound $\|u(t)\|_{H^s} \le c |t|^{\frac{2}{3}(s-1)+}$ $(s >> 1)$ for global solutions by using sharp bilinear estimates, which also imply new local-wellposedness results for rough data. Moreover in the case of blow-up solutions upper and lower bounds for $\|u(t)\|_{H^s}$ $(s>>1)$ were given.

We also consider the Cauchy problem for the Zakharov system in space dimension $n=2$ and $n=3$ :
\begin{eqnarray}
i\psi_t + \Delta \psi & = & \phi \psi \nonumber \\
\phi_{tt} - \Delta \phi & = & \Delta (|\psi|^2) \label{3.1}\\
\psi(0) \; = \; \psi_0 \quad , \quad \phi(0) & = & \phi_0 \quad , \quad \phi_t(0)\; = \;\phi_1 \nonumber
\end{eqnarray}
A general local existence and uniqueness result was proven by \cite{GTV} improving local (and global) wellposedness in energy space \cite{BC},\cite{C}.

Assuming $\psi_0 \in H^{s}$ , $ \phi_0 \in H^{m}$ and $ (-\Delta)^{-1/2}\phi_1 \in H^{m} $ with $s\ge 1$ , $ m > 0,$  $m+1>s>m$ it is not difficult to see that their arguments imply that the life span depends only on $\|\psi_0\|_{H^{1}}$ , $ \|\phi_0\|_{L^2} $ and $ \|(-\Delta)^{-1/2}\phi_1\|_{L^2}$ (cf. Prop. \ref{Proposition 3.1}). This solution can be extended globally by using the mass and energy conservation (cf.(\ref{3.1a}),(\ref{3.1b}),(\ref{3.1c})), provided $\|\psi_0\|_{L^2}$ is sufficiently small if $n=2$, and provided $\|\psi_0\|_{H^{1}} + \|\phi_0\|_{L^2} + \| (-\Delta)^{-1/2} \phi_1\|_{L^2} $ is sufficiently small if $n=3$. Again using the Fourier restriction norm method we show that in this situation a polynomial bound for $ \|\psi(t)\|_{H^{s}}$ $+ \|\phi(t)\|_{H^{m}}$ $+\|(-\Delta)^{-1/2}\phi_t(t)\|_{H^{m}}$ can be given for $ m+1>s>m $ and $ s >> 1 $ (Theorem \ref{Theorem 3.1}). This bound depends on $s$ and $m$. 

J. Colliander and G. Staffilani \cite{CS} considered the case $n=2$ for the Zakharov system and proved the bound $ \|\psi(t)\|_{H^s} \le c(1+t)^{s-1+} $ for $ s >> 1 $ and data in the Schwarz class and $\phi_1 \in \dot{H}^{-1} $ (\cite{CS}, Theorem 1). Their arguments imply that in fact the Schwarz space here can be replaced by the assumption $(\psi_0, \phi_0,\phi_1) \in H^{s} \times H^{s} \times H^{s-1}$ , and also a bound for $\|\phi(t)\|_{H^{s-1}} + \|\phi_t(t)\|_{H^{s-2}}$ can be given.

The necessary bilinear estimates for both problems are given in section 1, the KGS system is considered in section 2 and the Zakharov system in section 3.

For an equation of the form $iu_t - \varphi(-i\nabla_x)u = 0$ and $\varphi$ measurable we use the spaces $X^{s,b}_{\varphi}$ which are the completion of the Schwarz space ${\cal S}$ with respect to
\begin{eqnarray*}
\| f \|_{X^{s,b}_{\varphi}} : & = & \| <\xi>^s <\tau>^b {\cal F}(e^{-it\varphi(-i\nabla_x)}f(x,t)) \|_{L^2_{\xi\tau}} \\
& = & \| <\xi>^s <\tau+\varphi(\xi)>^b \widehat{f}(\xi,\tau) \|_{L^2_{\xi,\tau}}
\end{eqnarray*}
For $ \varphi(\xi) = \pm < \xi > $ or $ \pm |\xi| $ we use the notation $ X^{s,b}_{\pm} $ and for $ \varphi(\xi) = | \xi |^2 $ simply $ X^{s,b} $. For a given time interval $I$ we define 
$$ \| f \|_{X^{s,b}(I)} = \inf_{\widetilde{f}_{|I}=f} \| \widetilde{f} \|_{X^{s,b}} \quad \mbox{and similarly} \quad  \| f \|_{X^{s,b}_{\pm}(I)} $$ 
We rely on \cite{GTV} for the framework of the technique. For the method cf. also Bourgain \cite{B1},\cite{B2}, Klainerman and Machedon \cite{KM1},\cite{KM2}, and Kenig, Ponce and Vega \cite{KPV}.

Fundamental for our bilinear estimates are the following Strichartz type estimates for the Schr\"odinger equation:
$$ \| e^{it\Delta}\psi_0\|_{L^q_t(I,L^r_x({\bf R^n}))} \le c \|\psi_0\|_{L^2_x({\bf R^n})} $$
and
$$ \|f\|_{L^q_t(I,L^r_x({\bf R^n}))} \le c \|f\|_{X^{0,\frac{1}{2}+}(I)} $$
if $ 0 \le \frac{2}{q} = n(\frac{1}{2}-\frac{1}{r}) < 1 $ (cf. \cite{GTV}, Lemma 2.4). The endpoint case $(q,r)=(2, \frac{2n}{n-2}) $ is also admissible if $n \ge 3 $ (cf. \cite{KT}, Cor. 1.4).

We use the notation $<\lambda>:= (1+\lambda^2)^{1/2}$ for $ \lambda \in {\bf R} $ and $ a+ $ (resp. $ a- $) 
for a number slightly larger (resp. smaller) than $a$. Finally, $J^s := {\cal F}^{-1}<\xi>^s {\cal F} $ , where ${\cal F}$ denotes the spatial Fourier transform.\\
{\bf Acknowledgment:} The authors would like to thank J. Colliander and G. Staffilani for helpful correspondence concerning their work \cite{CS} on the Zakharov system.

\section{Bilinear estimates}
In this section we collect the estimates which are used to control the nonlinear terms in KGS as well as the Zakharov system.

We start with the estimates from \cite{GTV}, Lemma 3.4 and 3.5.
\begin{lemma}
\label{Lemma 1.1}
In space dimension $ 2 \le n \le 4 $ we have with $ \delta = 0$ for $n=2$ or $n=3$ and $\delta > 0$ arbitrarily small for $n=4$ :\\
a) $$ \|\psi_1 \overline{\psi_2} \|_{X_{\pm}^{m+1,-\frac{1}{2}+}} \le c \|\psi_1\|_{X^{s,\frac{1}{2}+}} \|\psi_2\|_{X^{s,\frac{1}{2}+}} $$
if $ s\ge 0$ , $ 2s \ge m+1+\delta $ , $ s>m $ .\\
b) $$ \|\phi_{\pm} \psi \|_{X^{s,-\frac{1}{2}+}} \le c \|\phi_{\pm}\|_{X_{\pm}^{m,\frac{1}{2}+}} \|\psi\|_{X^{s,\frac{1}{2}+}} $$
if $ m \ge \delta $ , $ 0 \le  s < m+1 $ .
\end{lemma}
Using the Leibniz rule for fractional derivatives we get as an immediate consequence the following 
\begin{lemma}
\label{Lemma 1.2}
In space dimension $ 2 \le n \le 4 $ the following estimates hold:\\  
a) if $ m \ge 3 $ and $\epsilon > 0$ arbitrarily small:
$$ \|\psi_1 \overline{\psi_2}\|_{X_{\pm}^{m-1,-\frac{1}{2}+}} \le c(\|\psi_1\|_{X^{m-2+2\epsilon,\frac{1}{2}+}} \|\psi_2\|_{X^{1-\epsilon,\frac{1}{2}+}} + \|\psi_1\|_{X^{1-\epsilon,\frac{1}{2}+}} \|\psi_2\|_{X^{m-2+2\epsilon,\frac{1}{2}+}}) $$
b) if $ s \ge 1 $ : 
$$ \|\phi_{\pm} \psi\|_{X^{s,-\frac{1}{2}+}} \le c(\| \phi_{\pm}\|_{X_{\pm}^{\delta,\frac{1}{2}+}} \| \psi\|_{X^{s,\frac{1}{2}+}} + \| \phi_{\pm}\|_{X_{\pm}^{s-1+,\frac{1}{2}+}} \| \psi\|_{X^{1-,\frac{1}{2}+}}) $$
with $ \delta \ge 0 $ , and $ \delta > 0 $ arbitrarily small in the case $n=4$. 
\end{lemma}

These estimates are sufficient to give a local existence result in the form we need it for our investigations, but in order to prove the desired polynomial growth for certain Sobolev norms of the solution we need another bilinear estimate which will be deduced next.
\begin{lemma}
\label{Lemma 1}
Let $ b > 1/2 $ . Then for space dimension $ 2 \le n \le 4 $ the following estimates hold true:
\begin{eqnarray}
\label{1}
\| fg \|_{X^{0,-b}} & \le & c \|f\|_{X^{0,b}} \|g\|_{X^{0,b}_{\pm}} \\
\label{2}
\| fg \|_{L_t^2(H_x^{-1-})} & \le & c \|f\|_{X^{0,b}} \|g\|_{X^{0,b}_{\pm}} \\
 \label{3}
\| fg \|_{X^{-1-,-b}} & \le & c \|f\|_{L^2_{xt}} \|g\|_{X^{0,b}_{\pm}} \\
\label{4}
\| fg \|_{X^{-1-,-b}} & \le & c \|f\|_{X^{0,b}} \|g\|_{L^2_{xt}}  
\end{eqnarray}
\end{lemma}
{\bf Proof:} Fix $ \frac{1}{q} = \frac{n}{8} $ , $ \frac{1}{q'} = 1-\frac{1}{q} $ , $ \frac{1}{p} = 1 - \frac{n}{4} $ , $ \frac{1}{r} = \frac{1}{2} + \frac{1}{n} $ and $ \frac{1}{\tilde{q}} = \frac{n}{4}-\frac{1}{2} $. From Strichartz' estimates we obtain the embeddings
\begin{equation}
\label{5}
X^{0,b} \subset L^q_t(L^4_x) 
\end{equation}
and
\begin{equation}
\label{6}
X^{0,b} \subset L^{\tilde{q}}_t(L^n_x)
\end{equation}
For $ n=4 $ this follows from the endpoint case $ (q,r) = (2,4) $ of the Strichartz estimates. Now using the dual version of (\ref{5}), H\"older, again (\ref{5}) as well as an embedding in the time variable, we obtain the following chain of inequalities:
$$
\|fg\|_{X^{0,-b}} \le c \|fg\|_{L^{q'}_t(L^{4/3}_x)} \le c \|f\|_{L_t^q(L^4_x)} \|g\|_{L^p_t(L^2_x)} \le c \|f\|_{X^{0,b}} \|g\|_{X^{0,b}_{\pm}} $$
This is (\ref{1}). To see (\ref{2}) we use Sobolev's embedding theorem in the $x$-variable, H\"older, (\ref{6}) and again a time embedding:
$$
\|fg\|_{L^2_t(H^{-1-})} \le c \|fg\|_{L^{2}_t(L^{r}_x)} \le c \|f\|_{L_t^{\tilde{q}}(L^n_x)} \|g\|_{L^p_t(L^2_x)} \le c \|f\|_{X^{0,b}} \|g\|_{X^{0,b}_{\pm}} $$
Concerning (\ref{3}) and (\ref{4}) we start by using the dual version of (\ref{5}) and Sobolev in $x$:
$$ \| fg \|_{X^{-1-,-b}} \le c \|fg\|_{L^{q'}_t(H^{-1-,\frac{4}{3}})} \le c \|fg\|_{L^{q'}_t(L^1_x)} $$
The latter is bounded by
$$ c \|f\|_{L^2_{xt}} \|g\|_{L^{2p}_t(L^2_x)} \le c \|f\|_{L^2_{xt}} \|g\|_{X^{0,b}_{\pm}} $$ 
which gives (\ref{3}), as well as by
$$ c \|f\|_{L^{2p}_t(L^2_x)} \|g\|_{L^2_{xt}} \le c \|f\|_{X^{0,b}} \|g\|_{L^2_{xt}} $$
leading to (\ref{4}).
\begin{lemma}
\label{Lemma 2}
Let $ b > 1/2 $ and $ 0 \le \Theta \le 1 $ . Then the estimate
$$ \|uv\|_{X^{-1,-b}} \le c \|u\|_{X^{-\Theta,b}} \|v\|_{X^{\Theta -1,b}_{\pm}} $$
holds, provided the space dimension $n$ fulfills $2\le n \le 4$ .
\end{lemma}
{\bf Proof:} We have
$$ \|uv\|_{X^{-1,-b}} = \| <\tau + |\xi|^2>^{-b} <\xi>^{-1} \int \widehat{u}(\xi_1,\tau_1) \widehat{v}(\xi_2,\tau_2) \, d\nu \|_{L^2_{\xi \tau}} $$
where $ d\nu = d\xi_1 \, d\tau_1 $ and $(\xi,\tau) = (\xi_1+\xi_2,\tau_1+\tau_2)$. 
We split the domain of integration into the regions $A$ and $B$, where in $A$ we have $|\xi_1| \le c |\xi| $ implying $|\xi_2| \le c |\xi| $. So the contribution from this part is bounded by
\begin{eqnarray*}
& & c \| <\tau + |\xi|^2>^{-b} \int <\xi_1>^{-\Theta} \widehat{u}(\xi_1,\tau_1) <\xi_2>^{\Theta -1} \widehat{v}(\xi_2,\tau_2) \, d\nu \|_{L^2_{\xi \tau}} \\
& & = \quad c \| (J^{-\Theta} u)(J^{\Theta -1}v) \|_{X^{0,-b}} \quad \le \quad c \|u\|_{X^{-\Theta,b}} \|v\|_{X^{\Theta -1,b}_{\pm}}
\end{eqnarray*}
by (\ref{1}). Next we consider the region $B$, where $|\xi| << |\xi_1| \sim |\xi_2| $ . This implies (for $\epsilon > 0$ sufficiently small)
$$|\xi|^{\epsilon} |\xi_1|^{\Theta} |\xi_2|^{1-\Theta} \le c ( <\tau + |\xi|^2>^b + <\tau_1+|\xi_1|^2>^b + <\tau_2 \pm |\xi_2| >^b ) $$
Thus we get three contributions from region $B$, the first of them being
$$ \| (J^{-\Theta} u)(J^{\Theta -1}v)\|_{L^2_t(H^{-1-\epsilon})} \le c \|u\|_{X^{-\Theta,b}} \|v\|_{X^{\Theta -1,b}_{\pm}} $$
by (\ref{2}). Writing $\Lambda^b = {\cal F}^{-1} <\tau + |\xi|^2>^b {\cal F}$ the second contribution 
(corresponding to the symbol $<\tau_1 + |\xi_1|^2>^b $) is
\begin{eqnarray*}
 \| (\Lambda^b J^{-\Theta} u)(J^{\Theta -1} v)\|_{X^{-1-,-b}} & \le & c \|\Lambda^b J^{-\Theta} u\|_{L^2_{xt}} \|J^{\Theta-1} v\|_{X^{0,b}_{\pm}} \\ & = & c \|u\|_{X^{-\Theta,b}} \|v\|_{X^{\Theta -1,b}_{\pm}} 
\end{eqnarray*}
where we used (\ref{3}). Finally, a similar argument using (\ref{4}) gives the same bound for the third contribution (corresponding to $<\tau_2 \pm |\xi_2|>^b$). 
\section{The Klein-Gordon-Schr\"odinger system}
We consider the Cauchy problem for the KGS system (\ref{2.1}) in space dimension $2 \le n \le 4$ . This system satisfies the conservation laws
\begin{equation}
\label{2.1a}
\|\psi(t)\|_{L^2} \equiv M
\end{equation}
and
\begin{equation}
\label{2.1b}
\|\nabla \psi(t)\|_{L^2}^2 + \frac{1}{2}(\|A^{1/2} \phi(t)\|_{L^2}^2 + \|\phi_t(t)\|_{L^2}^2) - \int |\psi(t)|^2 \phi(t) \, dx \equiv E
\end{equation}
where $ A := -\Delta + 1 $ .\\
By Gagliardo-Nirenberg we have for $6\ge n \ge 3$ :
$$ \left| \int_{{\bf R^n}} |\psi|^2 \phi \, dx \right| \le \|\phi\|_{L^{\frac{2n}{n-2}}} \|\psi\|_{L^{\frac{2n}{n-2}}}^{n-2} \|\psi\|_{L^2}^{6-n} \le c \|A^{1/2}\phi\|_{L^2} \|\nabla \psi\|_{L^2}^{n-2} \|\psi\|_{L^2}^{6-n} $$
If $n=3$ this is easily estimated by $\frac{1}{4}\|A^{1/2}\phi\|_{L^2}^2 + \frac{1}{2} \|\nabla \psi\|_{L^2}^2 +c \|\psi\|_{L^2}^6 $, if $n=4$ we have a bound by $\frac{1}{4} \|A^{1/2}\psi\|_{L^2}^2 + \frac{1}{2} \|\nabla \psi\|_{L^2}^2 $ , provided $\|\psi_0\|_{L^2}$ is sufficiently small. Because the case $n=2$ is easy to handle we have by the conservation laws in the case $2\le n \le 4 $ an a-priori bound
$$ \|\psi(t)\|_{H^{1}} + \|\phi(t)\|_{H^{1}} + \|\phi_t(t)\|_{L^2} \le const \quad \forall t \in {\bf R}^+ $$
The equivalent first order system reads as follows:
\begin{eqnarray}
\nonumber i\psi_t + \Delta \psi & = & - \frac{1}{2} (\phi_+ + \phi_-)\psi \\ 
\label{2.2} i\phi_{\pm t} \mp A^{1/2} \phi_{\pm} & = & \mp A^{-1/2}(|\psi|^2) \\
\nonumber \psi(0) \; = \; \psi_0 \quad , \quad \phi_{\pm}(0) & = & \phi_0 \pm iA^{-1/2} \phi_1 \; =: \; \phi_{0\pm}
\end{eqnarray}
where $ \phi_{\pm} = \phi \pm iA^{-1/2} \phi_t $ or, conversely, $\phi = \frac{1}{2}(\phi_+ + \phi_-) $ , $ 2iA^{-1/2}\phi_t = \phi_+ - \phi_- .$ \\
Using the corresponding system of integral equations 
\begin{eqnarray}
\nonumber
\psi(t) & = & e^{it\Delta} \psi_0 + i \int_0^t e^{i(t-s)\Delta} \frac{1}{2}(\phi_+(s) + \phi_-(s)) \psi(s) \, ds \\
\label{2.3}
\phi_{\pm}(t) & = & e^{\mp itA^{1/2}}\phi_{0\pm} \pm i \int_0^t e^{\mp i(t-s)A^{1/2}} A^{-1/2} (|\psi(s)|^2) \, ds
\end{eqnarray}
and Lemma \ref{Lemma 1.2} combined with the energy bound above we get the following version of the existence theorem by standard arguments:
\begin{prop}
\label{Proposition 2.1}
Let $ s \ge 1 $ , $ m > 0 $ satisfy $ s-1<m<s+2$ . Assume $\psi_0 \in H^{s} $ , $\phi_{0\pm} \in H^{m} $. Then there exists $T=T(\|\psi_0\|_{H^1},\|\phi_{0\pm}\|_{H^{\delta}}) > 0$ , where $ \delta = 0 $ , if $ 2\le n \le 3$ , and $\delta > 0 $ arbitrarily small, if $n=4$ , such that system (\ref{2.2}) has a unique solution $ \psi \in X^{s,\frac{1}{2}+}$ $[0,T] $ , $\phi_{\pm} \in X_{\pm}^{m,\frac{1}{2}+}[0,T] $ (especially $\psi \in C^0([0,T],H^{s}) $ ,$\phi_{\pm} \in C^0([0,T],H^{m}$)).
This solution satisfies:
$$
\|\psi\|_{X^{\sigma,\frac{1}{2}+}[0,T]} + \|\phi_{+}\|_{X_{+}^{\rho,\frac{1}{2}+}[0,T]} + \|\phi_{-}\|_{X_{-}^{\rho,\frac{1}{2}+}[0,T]} \le  d (\|\psi_0\|_{H^{\sigma}} + \|\phi_{0+}\|_{H^{\rho}} + \|\phi_{0-}\|_{H^{\rho}}) 
$$
for $ \sigma \ge 1 $ , $ \rho > 0 $ and $ \sigma -1 < \rho < \sigma +2 $ . 
This solution exists globally, provided $\|\psi_0\|_{L^2}$ is sufficiently small in dimension $n=4$.
\end{prop} 
{\bf Remark:} A global existence result for rougher data also holds true, namely $ (\psi_0,\phi_{0\pm}) \in H^{s} \times H^{m} $ with $ 1 \ge s,m > 7/10 $ and $ s+m > 3/2 $ (cf. \cite{P}).

Our aim is to give a bound in time for $ \|\psi(t)\|_{H^{s}}$ and $\|\phi_{\pm}(t)\|_{H^{m}}$ for $s,m >>1$.

The proof uses the following well-known elementary observation (used by \cite{B},e.g.).
\begin{lemma}
\label{Lemma 2.2}
Assume $f \in C^0({\bf R^+},{\bf R^+})$. Moreover assume the existence of $t_0 > 0$ s. th. $\forall \tau,t \in [\tau,\tau + t_0] $ the following estimate holds
$$ f(t) \le f(\tau) + c f(\tau)^{1-\delta} $$
where $ 0 < \delta \le 1 $ , $ c > 0 $ . Then there exist $c_1,c_2 > 0$ s. th. 
$$ f(t) \le c_1 f(0) + c_2 t^{\frac{1}{\delta}} $$
for $ t \ge 0 $ .
\end{lemma}
Our main result for the KGS system is the following:
\begin{theorem}
\label{Theorem 2.1}
Let $ s > 1 $ be an even integer, $ m \ge 3 $ and $ s+2 > m > s-1 $ . For $ 2 \le n \le 4 $ assume that $(\psi_0,\phi_0,\phi_1) \in H^{s}\times H^{m} \times H^{m-1} $ with $\|\psi_0\|_{L^2}$ sufficiently small in the case $n=4$ . Then the global solution of Prop. \ref{Proposition 2.1} for the KGS system (\ref{2.1}) fulfills:
\begin{equation}
\|\psi(t)\|_{H^{s}} + \|\phi(t)\|_{H^{m}} + \|\phi_t(t)\|_{H^{m-1}}  \le c (1+t)^{\frac{1}{\delta}}   
\label{2.4}
\end{equation}
where $ \frac{1}{\delta} = \max(\frac{m-1}{m-s+1},s-1,\frac{s-1}{s-m+2})$ . If in addition $ s \le m \le s+1 $ , we have $ \frac{1}{\delta} = s-1 $ .
\end{theorem} 
{\bf Proof:} We use the KGS system (\ref{2.2}) on any interval $I$ where the local existence theorem Prop. \ref{Proposition 2.1} applies. We assume w.l.o.g. $I=[0,T].$ \\
{\bf Estimate for the Schr\"odinger part:}
\begin{eqnarray*}
\|\psi(t)\|_{H^s}^2 - \|\psi_0\|_{H^s}^2 = \int_0^t \frac{\partial}{\partial \tau} \|\psi(\tau)\|_{H^s}^2 \, d\tau = 2 Re \int_0^t <\psi_t(\tau),\psi(\tau)>_{H^s} \, d\tau = \\
2 Re (i  \int_0^t \hspace{-0.5em}<\Delta \psi(\tau),\psi(\tau)>_{H^s}  d\tau ) + 2 Re (i \int_0^t \hspace{-0.5em} <(1-\Delta)^{\frac{s}{2}}(\psi \phi),(1-\Delta)^{\frac{s}{2}}\psi>_{L^2}  d\tau )
\end{eqnarray*}
The first term vanishes, since $<\Delta \psi(\tau),\psi(\tau)>_{H^s} \in {\bf R} $ , and the second one can be estimated by
$$ c \sum_{|\alpha_1|+|\alpha_2|\le s \, , \, |\alpha_2| < s} \int_I |<D^{\alpha_1}\phi(\tau) D^{\alpha_2}\psi(\tau),(1-\Delta)^{\frac{s}{2}} \psi(\tau)>_{L^2_x}| \, d\tau $$
In this sum no $|\alpha_2|=s$ - contribution occurs, since $<\phi(\tau),|(1-\Delta)^{\frac{s}{2}} \psi|^2>_{L^2_x} \in {\bf R}.$  The $|\alpha_1|=s$ - contributions can be estimated by Lemma \ref{Lemma 2} with $\theta = 0$ and Prop. \ref{Proposition 2.1}:
\begin{eqnarray} \nonumber
\lefteqn{\hspace{-2em} \int_I |<\psi(\tau),D^{\alpha_1}\phi(\tau) (1-\Delta)^{\frac{s}{2}}\psi(\tau)>_{L^2_x}| \, d\tau  \le  \|\psi\|_{X^{1,b}(I)} \|J^s \phi J^s \psi\|_{X^{-1,-b}(I)}} \\ \label{****}
& \le & c \|\psi\|_{X^{1,b}(I)} (\|J^{s-1} \phi_+\|_{X^{0,b}_+(I)} + \| J^{s-1} \phi_-\|_{X^{0,b}_-(I)}) \|\psi\|_{X^{s,b}(I)} \\ \nonumber
& \le & c (\|\psi_0\|_{H^1} + \|(\phi_{0+},\phi_{0-})\|_{H^1})(\|(\phi_{0+},\phi_{0-})\|_{H^{s-1}} + \|\psi_0\|_{H^{s-3+}}) \\ \nonumber
& & \times  (\|\psi_0\|_{H^s} + \|(\phi_{0+},\phi_{0-})\|_{H^{s-1+}})
\end{eqnarray}
(here $\|(\phi_{0+},\phi_{0-})\|_{H^{s-1}}^2 := \|\phi_{0+}\|_{H^{s-1}}^2 + \|\phi_{0-}\|_{H^{s-1}}^2 $ ). Proceeding similarly with (\ref{1}) instead of Lemma \ref{Lemma 2} we get the following bound for the remaining terms:
\begin{eqnarray}
\nonumber 
& \hspace{-1em} c & \sum_{|\alpha_1| + |\alpha_2| \le s \, , \, |\alpha_1|,|\alpha_2| < s} (\|\phi_+\|_{X^{|\alpha_1|,b}_+(I)} + \|\phi_-\|_{X^{|\alpha_1|,b}_-(I)})\|\psi\|_{X^{|\alpha_2|,b}(I)} \|\psi\|_{X^{s,b}(I)} \\ \label{***}
& \hspace{-1em}  \le &  \hspace{-1em} c \|\psi\|_{X^{s,b}(I)}((\|\phi_+\|_{X_+^{s-1,b}(I)} + \|\phi_-\|_{X_-^{s-1,b}(I)}) \|\psi\|_{X^{1,b}(I)} \\ \nonumber & & \qquad \qquad +(\|\phi_+\|_{X_+^{1,b}(I)} + \|\phi_-\|_{X_-^{1,b}(I)}) \|\psi\|_{X^{s-1,b}(I)})  \\ \nonumber
& \hspace{-1em} \le &  \hspace{-1em} c(\|\psi_0\|_{H^s} + \|(\phi_{0+},\phi_{0-})\|_{H^{s-1+}})[(\|(\phi_{0+},\phi_{0-})\|_{H^{s-1}} + \|\psi_0\|_{H^{s-3+}})(\|\psi_0\|_{H^1} \\ \nonumber & \hspace{-0.8em} + &  \hspace{-0.8em} \|(\phi_{0+},\phi_{0-})\|_{H^1})+
 (\|(\psi_{0+},\phi_{0-})\|_{H^1} \hspace{-0.2em} + \hspace{-0.2em} \|\psi_0\|_{H^{1}})(\|\psi_0\|_{H^{s-1}} + \|(\phi_{0+},\phi_{0-})\|_{H^{s-2+}})]
\end{eqnarray}
Taking into account that $\|\psi_0\|_{H^1} + \|(\phi_{0+},\phi_{0-})\|_{H^1} \le const $ we arrive at
\begin{eqnarray}
\label{*}
\lefteqn{\|\psi(t)\|_{H^s}^2 - \|\psi_0\|_{H^s}^2} \\
\nonumber & \le & c(\|(\phi_{0+},\phi_{0-})\|_{H^{s-1}} + \|\psi_0\|_{H^{s-1}})(\|\psi_0\|_{H^s} + \|(\phi_{0+},\phi_{0-})\|_{H^m}) \\
\nonumber
& \le & c(\|(\phi_{0+},\phi_{0-})\|_{H^m}^{1-\vartheta_1} + \|\psi_0\|_{H^s}^{1-\vartheta_2})(\|\psi_0\|_{H^s} + \|(\phi_{0+},\phi_{0-})\|_{H^m})
\end{eqnarray}
with $\vartheta_1 = \frac{m-s+1}{m-1} $ ( $ \in (0,1] $ , provided $ m>s-1\ge 1$) and with $\vartheta_2 = \frac{1}{s-1} $ .\\
{\bf Estimate for the Klein-Gordon part:} a similar computation as above shows that
\begin{eqnarray} \nonumber
\lefteqn{\|(\phi_+(t),\phi_-(t))\|_{H^m}^2 - \|(\phi_{0+},\phi_{0-})\|_{H^m}^2} \\
\nonumber  & = & 2 Re (i\int_0^t <J^{m-1}(|\psi(\tau)|^2),J^m(\phi_+(\tau)-\phi_-(\tau))>_{L^2_x} d\tau ) \\ \nonumber
& \le & c \|\psi\|_{X^{1,b}(I)} (\|J^{m-1}\psi \, J^m \phi_+\|_{X^{-1,-b}(I)} + \|J^{m-1}\psi \, J^m \phi_-\|_{X^{-1,-b}(I)}) \\ \nonumber
& \le & c \|\psi\|_{X^{m-2,b}(I)} (\|\phi_+\|_{X_+^{m,b}(I)} + \|\phi_-\|_{X_-^{m,b}(I)}) \\ \nonumber
& \le & c (\|\psi_0\|_{H^{m-2}}+\|(\phi_{0+},\phi_{0-})\|_{H^{m-3+}}) (\|(\phi_{0+},\phi_{0-})\|_{H^m} + \|\psi_0\|_{H^{m-2+}}) \\ \label{**}
& \le & c(\|\psi_0\|_{H^s}^{1-\vartheta_3} + \|(\phi_{0+},\phi_{0-})\|_{H^m}^{1-\vartheta_4})(\|(\phi_{0+},\phi_{0-})\|_{H^m} + \|\psi_0\|_{H^s})
\end{eqnarray}
with $\vartheta_3 = \frac{s-m+2}{s-1}$ and $\vartheta_4 = \frac{3}{m-1}+$ ,
where we now used Lemma \ref{Lemma 2} with $\Theta =1$. \\
{\bf Conclusion:} combining (\ref{*}) and (\ref{**}) we see that for
$$ f(t) = \|\psi(t)\|_{H^s}^2 + \|(\phi_+(t),\phi_-(t))\|_{H^m}^2 $$
the inequality
$$ f(t) \le f(0) + c f(0)^{1-\frac{\vartheta}{2}} $$ 
holds with $ \vartheta = \min_{1 \le i \le 3} \vartheta_i $ , provided $s+2>m\ge 3$ and $m>s-1\ge 1$ (remark that $\vartheta_1 < \vartheta_4$). Lemma \ref{Lemma 2.2} gives $ f(t) \le c (1+t)^{\frac{2}{\delta}} $ , where $ \frac{1}{\delta} = \max(\frac{m-1}{m-s+1},s-1,\frac{s-1}{s-m+2}) $ . This implies the claimed result.

\section{The Zakharov system}
The Cauchy problem for the Zakharov system (\ref{3.1}) is considered in space dimension $n=2$ and $n=3$.\\
The conservation laws are the following:
\begin{equation}
\label{3.1a}
\|\psi(t)\|_{L^2} \equiv M
\end{equation}
and
\begin{equation}
\label{3.1b}
E(\psi,\phi) \equiv \int (|\nabla \psi|^2 + \frac{1}{2}(|\phi|^2+|V|^2)+\phi|\psi|^2) \, dx
\end{equation}
where $ div \, V = \phi_t $ .\\
It is well-known (cf. \cite{SS}, Proof of Thme. 1) that this implies a uniform a-priori-bound
\begin{equation}
\label{3.1c}
\|\psi\|_{H^{1}} + \|\phi\|_{L^2} + \| (-\Delta)^{-1/2} \phi_t \|_{L^2} \le K \qquad \forall t \in {\bf R^+}
\end{equation}
if $ \|\psi_0\|_{L^2} $ is sufficiently small for $n=2$, and if $ \|\psi_0\|_{H^{1}} + \|\phi_0\|_{L^2} + \|(-\Delta)^{-1/2} \phi_1\|_{L^2} $ is sufficiently small for $n=3$ . \\
The equivalent first order system is
\begin{eqnarray}
\nonumber
i \psi_t + \Delta \psi & = & \frac{1}{2} (\phi_+ + \phi_-) \psi \\
\label{3.2}
i \phi_{\pm t} \mp (-\Delta)^{1/2} \phi_{\pm} & = & \pm (-\Delta)^{1/2} (|\psi|^2) \\
\nonumber
\psi(0) \; = \; \psi_0 \quad , \quad \phi_{\pm}(0) & = & \phi_0 \pm i (-\Delta)^{-1/2} \phi_1 \; =: \; \phi_{0\pm} 
\end{eqnarray}
where $ \phi_{\pm} = \phi \pm i(-\Delta)^{-1/2}\phi_t $ or, conversely, $\phi = \frac{1}{2}(\phi_+ + \phi_-) $ , $ 2i(-\Delta)^{-1/2} \phi_t = \phi_+ - \phi_- $ .\\
The corresponding system of integral equations is
\begin{eqnarray}
\nonumber
\psi(t) & = & e^{it\Delta} \psi_0 - i \int_0^t e^{i(t-s)\Delta} \frac{1}{2}(\phi_+(s) + \phi_-(s)) \psi(s) \, ds \\
\label{3.3}
\phi_{\pm}(t) & = & e^{\mp it(-\Delta)^{1/2}}\phi_{0\pm} \mp i \int_0^t e^{\mp i(t-s)(-\Delta)^{1/2}} (-\Delta)^{1/2} (|\psi(s)|^2) \, ds
\end{eqnarray}
Lemma \ref{Lemma 1.1} and \ref{Lemma 1.2} combined with the energy bound above implies the following existence theorem by standard arguments:
\begin{prop}
\label{Proposition 3.1}
Let $ n=2$ or $n=3$ , $ s \ge 1 $ , $ m > 0 $ satisfy $ m+1 > s > m $ . Assume $\psi_0 \in H^{s} $  , $ \phi_{0\pm} \in H^{m}$ . Then there exists $ T = T(\|\psi_0\|_{H^{1-}}, \|\phi_{0\pm}\|_{L^2}) > 0 $ such that the system (\ref{3.2}) has a unique solution $ \psi \in X^{s,\frac{1}{2}+}[0,T]$ , $ \phi_{\pm} \in X_{\pm}^{m,\frac{1}{2}+}[0,T] $ . This solution satisfies 
$$
\|\psi\|_{X^{\sigma,\frac{1}{2}+}[0,T]} + \|\phi_{+}\|_{X_{+}^{\rho,\frac{1}{2}+}[0,T]} + \|\phi_{-}\|_{X_{-}^{\rho,\frac{1}{2}+}[0,T]} \le  d (\|\psi_0\|_{H^{\sigma}} + \|\phi_{0+}\|_{H^{\rho}} + \|\phi_{0-}\|_{H^{\rho}})
$$
for $ 1 \le \sigma \le s $ , $ 0 < \rho \le m $ and $\sigma > \rho > \sigma -1$.
This solution exists globally, provided $ \|\psi_0\|_{L^2} $ is sufficiently small for $n=2$ , and provided $\|\psi_0\|_{H^{1}} + \|\phi_{0+}\|_{L^2} + \|\phi_{0-}\|_{L^2} $ is sufficiently small for $n=3$.
\end{prop}
{\bf Remark:} A more general local result for rougher data can also be given (cf. \cite{GTV}).

Our aim is to give a bound in time for $\|\psi(t)\|_{H^{s}} $ and $ \|\phi_{\pm}(t)\|_{H^{m}} $ for $ s,m >> 1 $. The main result is the following
\begin{theorem}
\label{Theorem 3.1}
Let $ n=2 $ or $ n=3 $ , let $ s > 0 $ (an even integer) and $ m \ge 0 $ satisfy $ m+1 > s > m $ . Assume $ \psi_0 \in H^{s} $  , $ \phi_0 \in H^{m} $ , $ (-\Delta)^{-1/2} \phi_1 \in H^{m}$ with $ \|\phi_0\|_{L^2} $ sufficiently small for $n=2$, and with $ \|\psi_0\|_{H^{1}} + \|\phi_0\|_{L^2} + \|(-\Delta)^{-1/2} \phi_1\|_{L^2} $ sufficiently small for $n=3$ . Then the global solution of Prop. \ref{Proposition 3.1} for the Zakharov system fulfills
$$ \|\psi(t)\|_{H^{s}} + \|\phi(t)\|_{H^{m}} + \| (-\Delta)^{-1/2} \phi_t (t)\|_{H^{m}} \le c (1+t)^{\frac{1}{\delta}} $$
where $ \frac{1}{\delta} = \max (\frac{m(s-1)}{m-s+1} , \frac{s-1}{s-m}+) $ .
\end{theorem}
{\bf Remark:} Especially, for $ m=s-\frac{1}{2}$ we have $\frac{1}{\delta} = (s-1)(2s-1) = $ {\cal O}$(2s^2)$ , and for $ m=s-\frac{1}{s} $ we have $\frac{1}{\delta} = $ {\cal O}$(s^2)$ .\\
{\bf Proof:} {\bf Estimate for the Schr\"odinger part:} In the proof of Theorem \ref{Theorem 2.1} we have shown in (\ref{****}) and (\ref{***}):
\begin{eqnarray*}
\|\psi(t)\|_{H^s}^2 - \|\psi_0\|_{H^s}^2  & \le & c (\|\psi\|_{X^{1,b}(I)} (\|\phi_+\|_{X_+^{s-1,b}(I)} + \|\phi_-\|_{X_-^{s-1,b}(I)}) \\ & &  + \|\psi\|_{X^{s-1,b}(I)}(\|\phi_+\|_{X_+^{1,b}(I)} + \|\phi_-\|_{X_-^{1,b}(I)})) \|\psi\|_{X^{s,b}(I)} 
\end{eqnarray*}
By use of Prop. \ref{Proposition 3.1} and (\ref{3.1c}) the latter is estimated by:
\begin{eqnarray}
\nonumber & & \hspace{-1.8em}
c[(\|\psi_0\|_{H^1} + \|(\phi_{0+},\phi_{0-})\|_{H^{0+}})((\|(\phi_{0+},\phi_{0-})\|_{H^{s-1}} + \|\psi_0\|_{H^{s-1+}} )) + (\|\psi_0\|_{H^{s-1}} \\ \nonumber  & & \hspace{-1.5em}+ \|(\phi_{0+},\phi_{0-})\|_{H^{s-2+}} (\|(\phi_{0+},\phi_{0-})\|_{H^1} + \|\psi_0\|_{H^{1+}})] (\|\psi_0\|_{H^s}+\|(\phi_{0+},\phi_{0-})\|_{H^{s-1+}}) \\
\nonumber
& &  \hspace{-1.8em} \le c[(1+\|(\phi_{0+},\phi_{0-})\|_{H^m}^{0+})(\|(\phi_{0+},\phi_{0-})\|_{H^m}^{\frac{s-1}{m}} + \|\psi_0\|_{H^s}^{\frac{s-2}{s-1}+}) \\ \nonumber
 & & +(\|\psi_0\|_{H^s}^{\frac{s-2}{s-1}} + \|(\phi_{0+},\phi_{0-})\|_{H^m}^{\frac{s-2}{m}+})(\|(\phi_{0+},\phi_{0-})\|_{H^m}^{\frac{1}{m}} + \|\psi_0\|_{H^s}^{0+})] \\
\nonumber
& & \qquad \qquad \times (\|\psi_0\|_{H^s} + \|(\phi_{0+},\phi_{0-})\|_{H^m}) \\ \nonumber
& &  \hspace{-1.8em} \le c(\|(\phi_{0+},\phi_{0-})\|_{H^m}^{\frac{s-1}{m}+} + \|\psi_0\|_{H^s}^{\frac{s-2}{s-1}+} + \|\psi_0\|_{H^s}^{\frac{s-2}{s-1}} \|(\phi_{0+},\phi_{0-})\|_{H^m}^{\frac{1}{m}})\\
\nonumber
& & \qquad \qquad \times (\|\psi_0\|_{H^s} + \|(\phi_{0+},\phi_{0-})\|_{H^m}) \\ \label{*****}
& &  \hspace{-1.8em} \le c(\|(\phi_{0+},\phi_{0-})\|_{H^m}^{1-\vartheta_1} + \|\psi_0\|_{H^s}^{1-\vartheta_1})(\|\psi_0\|_{H^s} + \|(\phi_{0+},\phi_{0-})\|_{H^m})
\end{eqnarray}
where $ \vartheta_1 = \frac{m-s+1}{m(s-1)} $ .\\
{\bf Estimate for the wave part:} In the estimate for the Klein-Gordon part in the proof of Theorem \ref{Theorem 2.1} we replace the term $J^{-1}(|\psi|^2)$ by $J^{-1}\Delta(|\psi|^2)$ and arrive at
\begin{eqnarray} \nonumber
\lefteqn{\|(\phi_+(t),\phi_-(t))\|_{H^m}^2 - \|(\phi_{0+},\phi_{0-})\|_{H^m}^2} \\
 \nonumber
& \le & c \|\psi\|_{X^{1,b}(I)} \|\psi\|_{X^{m,b}(I)} (\|\phi_+\|_{X_+^{m,b}(I)} + \|\phi_-\|_{X_-^{m,b}(I)}) \\ \nonumber
& \le & c (\|\psi_0\|_{H^1} + \|(\phi_{0+},\phi_{0-})\|_{H^{0+}})(\|\psi_0\|_{H^{m}}+\|(\phi_{0+},\phi_{0-})\|_{H^{m-1+}}) \\ & & \hspace{13em} \times (\|(\phi_{0+},\phi_{0-})\|_{H^m} + \|\psi_0\|_{H^{m+}}) \\ \label{******}
& \le & c(\|\psi_0\|_{H^s}^{1-\vartheta_2} + \|(\phi_{0+},\phi_{0-})\|_{H^m}^{1-\vartheta_3})(\|(\phi_{0+},\phi_{0-})\|_{H^m} + \|\psi_0\|_{H^s})
\end{eqnarray}
where we used (\ref{3.1c}) and Prop. \ref{Proposition 3.1}. Now $\vartheta_2 = \frac{s-m}{s-1}-$ and $\vartheta_3 = \frac{1}{m}-$ .
Adding (\ref{*****}) and (\ref{******}) we get:
\begin{eqnarray*}
\lefteqn{\|\psi(t)\|_{H^s}^2 - \|\psi_0\|_{H^s}^2 + \|(\phi_{0+},\phi_{0-})\|_{H^m}^2 - \|(\phi_{0+},\phi_{0-})\|_{H^m}^2} \\
& \le & c(\|(\phi_{0+},\phi_{0-})\|_{H^m}^{1-\vartheta_1} + \|\psi_0\|_{H^s}^{1-\vartheta_1} + \|\psi_0\|_{H^s}^{1-\vartheta_2} + \|(\phi_{0+},\phi_{0-})\|_{H^m}^{1-\vartheta_3})\\
& & (\|(\phi_{0+},\phi_{0-})\|_{H^m} + \|\psi_0\|_{H^s})
\end{eqnarray*}
Using Lemma \ref{Lemma 2.2} this implies the claimed estimate similarly as for the KGS system.

\end{document}